\documentclass[11pt,a4paper]{article}
\usepackage{style}%
\usepackage{amsthm}
\usepackage[capitalize]{cleveref}

\crefname{equation}{}{}
\Crefname{equation}{}{}
\crefname{figure}{Figure}{Figures}
\author{Stefano Berrone, Andrea Borio, Francesca Marcon, Gioana Teora \footnote{The authors are members of the INdAM-GNCS.  The authors kindly
    acknowledge partial financial support provided by INdAM-GNCS Projects 2022, by the
    MIUR project ``Dipartimenti di Eccellenza'' Programme (2018–2022)
    CUP:E11G18000350001 and by the PRIN 2020 project (No. 20204LN5N5\_003).}}%
\title{A first-order stabilization-free Virtual Element Method}%
\date{}

\begin{document}


\maketitle
\begin{abstract}
 In this paper, we introduce a new Virtual Element Method (VEM) not requiring any  stabilization term based on the usual enhanced first-order VEM space. The new method relies on a modified formulation of the discrete diffusion operator that ensures stability preserving all the properties of the differential operator.
\end{abstract}

\section{Introduction}

Recently, in the context of Virtual Element Methods (VEM), a growing interest
has been devoted to the definition of bilinear forms not
requiring a stabilization term. In \cite{BBME2VEM}, a lowest-order
stabilization-free scheme was proposed and analysed, proving that it is possible
to define coercive bilinear forms based on polynomial projections of virtual
basis functions of suitable high-degree polynomial spaces.
In \cite{BERRONE2022107971}, the
proposed scheme was compared to standard VEM, and results showed that the
absence of a stabilization operator can reduce the error and help convergence in
case of strongly anisotropic problems.

In this paper, we propose a variation of
the scheme introduced in \cite{BBME2VEM}, strongly exploiting the theory
developed in that paper to choose the smallest possible polynomial space that
guarantees coercivity.


We consider an open bounded domain $\Omega \subset \mathbb{R}^2$ and the
following standard advection-diffusion-reaction problem:
{\sl find $u\in\sobh[0]{1}{\Omega}$ such that}
\begin{equation}
  \label{eq:modelvarform}
  \scal[\Omega]{\K\nabla u}{\nabla v} + \scal[\Omega]{\bs\beta\cdot\nabla u }{v} +
  \scal[\Omega]{\gamma u }{v} = \scal[\Omega]{f}{v}
  \quad \forall v \in \sobh[0]{1}{\Omega} \,,
\end{equation}
where $\scal[\Omega]{\cdot}{\cdot}$ denotes the $\lebl{\Omega}$ scalar
product. We make standard assumptions on the coefficients in order to guarantee
the well-posedness of the problem, namely, all coefficients are $\lebl[\infty]{\Omega}$, $\K$ is a symmetric uniformly
positive definite tensor, $\div\bs\beta = 0$, and
$\inf_{x\in\Omega}\gamma(x) \geq 0$. Here we consider homogeneous Dirichlet
boundary conditions, 
but more general boundary conditions can be considered.

\section{Local spaces and projections}
\label{sec:local}
We consider a family of polygonal tessellations $\Mh$ of $\Omega$, satisfying
the following standard mesh assumptions: $\exists \kappa > 0$ such that $\forall E \in \Mh$, 
$E$ is star-shaped with respect to a ball of radius $\rho \geq \kappa h_E$, and $\forall e \in \E{E}$, where $\E{E}$ is the set of edges of $E$, $\abs{e} \geq \kappa h_E$,
where $h_E$ denotes the diameter of $E$.
For any given $E\in\Mh$, we define the following standard Virtual Element space
\cite{Beirao2015b}:
\begin{multline*}
  \Vh[E] = \left\{v\in\sobh{1}{E}\colon \Delta v \in \Poly{1}{E},\,
    \trace{v}{\partial E}{}\in \cont{\partial E}, \trace{v}{e}{} \in \Poly{1}{e}
    \; \forall e \in \E{E},\right.
  \\
  \left.\scal[E]{v-\proj[\nabla,E]{1}{}v}{p} = 0 \;\forall p \in \Poly{1}{E}
  \right\} \,,
\end{multline*}
where $\trace{v}{\omega}{}$ denotes the trace of $v$ on $\omega$ and
$\proj[\nabla,E]{1}{} v \in \Poly{1}{E}$ is defined such that
$ \scal[E]{\nabla v-\nabla \proj[\nabla,E]{1}{}v}{\nabla p} = 0$
$\forall p \in \Poly{1}{E}$ and
$\int_{\partial E}v = \int_{\partial E}\proj[\nabla,E]{1}{}v$.  The degrees of
freedom of $\Vh[E]$ are the values of functions at the vertices of the polygon
$E$. 


For any given $\ell \in \mathbb{N}$, we define the following spaces of
harmonic polynomials of degree $\ell + 1$:
\begin{equation*}
\HPoly{\ell+1}{E} = \left\{ p \in \Poly{\ell+1}{E}\colon \Delta p = 0, \int_Ep = 0 \right\}\,.
\end{equation*}
Let $\nabla\HPoly{\ell+1}{E}$ be the space of gradients of functions in
$\HPoly{\ell+1}{E}$. We define the projector
$\proj[\HPoly{}{},E]{\ell}{}\colon \lebldouble{E} \to \nabla \HPoly{\ell+1}{E}$
such that, $\forall \bs{v} \in \lebldouble{E}$,
\begin{equation}
  \label{eq:defPiNablaH}
  \scal[E]{\proj[\HPoly{}{},E]{\ell}{} \bs{v}}{\nabla p_{\ell+1}} =
  \scal[E]{\bs{v}}{\nabla p_{\ell+1}} \quad  \forall p_{\ell+1} \in \HPoly{\ell+1}{E}\,.
\end{equation}
Notice that, since $\HPoly{\ell+1}{E}$ does not contain constants by definition,
$\nabla p_{\ell+1}$ is never zero in \eqref{eq:defPiNablaH} and
$\dim\nabla\HPoly{\ell+1}{E} = \dim\HPoly{\ell+1}{E} = 2\ell + 2$. Moreover, notice
that $\PolyDouble{0}{E} \subseteq \nabla \HPoly{\ell+1}{E}$, and in particular
$\PolyDouble{0}{E} = \nabla \HPoly{1}{E}$.

Now, given a function $v_h \in \Vh[E]$, consider the problem of computing
$\proj[\HPoly{}{},E]{\ell}{}\nabla v_h$. 
Let $\left\{ h_i , i = 1,\ldots,2\ell+2\right\}$ be a set of basis
functions of $\HPoly{\ell+1}{E}$. Then
$\proj[\HPoly{}{},E]{\ell}{}\nabla v_h = \sum_{j=1}^{2l+2} d_j \nabla h_j$,
where the values $d_j$ can be computed by solving the following system of
equations:
\begin{equation}
  \label{eq:systemPiNablaH}
  \sum_{j=1}^{2l+2} \scal[E]{\nabla h_j}{\nabla h_i} d_j =
  \scal[E]{\nabla v_h}{\nabla h_i} \quad \forall i = 1,\ldots,2\ell+2 \,.
\end{equation}
The right-hand side can be computed since we know $v_h$ analitically on the boundary, recalling that $\Delta h_i = 0$ and applying Green's theorem:
$
  \scal[E]{\nabla v_h}{\nabla h_i} = \scal[\partial E]{v_h}{\pder{h_i}{n}},
  \
  \forall i = 1,\ldots,2\ell+2 \,.
  $
On each edge, the right-hand side is the integral of a polynomial of degree $\ell+1$, that can be computed exactly using $\lceil\frac{\ell+2}{2}\rceil$ Gauss quadrature nodes.
%
Concerning the left-hand side of \eqref{eq:systemPiNablaH}, a way to reduce the computational cost, with respect to 2D quadrature rules, 
is to observe that
$
  \scal[E]{\nabla h_j}{\nabla h_i} = \scal[\partial E]{h_j}{\pder{h_i}{n}} \,,
  $
that is the integral of a piecewise polynomial of degree $2\ell+1$. Then, the
integral can be computed by $\ell+1$ Gauss quadrature nodes on each edge,
reducing the number of function evaluations to $\sim N_E \ell$.

\section{Discrete variational formulation}
\label{sec:discrvarform}

Let
$\Vh = \{v_h \in \sobh[0]{1}{\Omega}\colon v_h \in \Vh[E] \;\forall E \in \Mh\}$
and let $\ell_E \geq 0$ be given $\forall E \in \Mh$, possibly different from
one polygon to another. Then, we look for $u_h \in \Vh$ such that
\begin{equation}
  \label{eq:discrvarform}
  \begin{split}
    \sum_{E\in\Mh}\scal[E]{\K\proj[\HPoly{}{},E]{\ell_E}{}\nabla u_h}
    {\proj[\HPoly{}{},E]{\ell_E}{}\nabla v_h} +
    \scal[E]{\bs\beta\cdot\proj[\HPoly{}{},E]{\ell_E}{}\nabla u_h } {
      \proj[0,E]{0}{}v_h}
    \\
    +\scal[E]{\gamma \proj[0,E]{0}{}u_h }{\proj[0,E]{0}{} v_h} = \sum_{E\in\Mh}
    \scal[E]{f}{\proj[0,E]{0}{} v_h} \quad \forall v_h \in \Vh \,,
  \end{split}
\end{equation}
where $\proj[0,E]{0}{}$ is the $\lebl{}$ projection operator onto constants. The following result provides the crucial ingredient for the
well-posedness of \eqref{eq:discrvarform}.
\begin{theorem}
  \label{th:wellposed}
  Assume that, $\forall E \in \Mh$, $2\ell_E + 2 \geq N_E - 1$, $N_E$ being the
  number of vertices of $E$. Then there exist $\alpha^\ast, \alpha_\ast > 0$,
depend on the mesh regularity parameter $\kappa$ and on local variations of $\K$,
  such that, $\forall u_h \in \Vh$, $\forall E \in \Mh$,
  \begin{equation*}
    \alpha_\ast\norm[E]{\sqrt{\K} \nabla u_h} \leq \norm[E]{\sqrt{\K}\proj[\HPoly{}{},E]{\ell_E}{}\nabla u_h} \leq \alpha^\ast\norm[E]{\sqrt{\K}\nabla u_h} \,.
  \end{equation*}
\end{theorem}
\begin{proof}
  The result follows from the theory developed in \cite{BBME2VEM}.
\end{proof}
Theorem \ref{th:wellposed} provides us a sufficient condition for the coercivity
of the diffusivity term of \eqref{eq:discrvarform}. The well-posedness of the
discrete problem is then obtained by the same arguments as in
\cite{Beirao2015b}. Optimal order a priori error estimates can be proved using
the techniques in \cite{Beirao2015b,BBME2VEM}. In particular, we get
\begin{align*}
  \norm[\Omega]{\sqrt{\K}\nabla (u-u_h)} &= O(h) \,,
  & \norm[\Omega]{u-u_h} = O(h^2) \,.
\end{align*}

\begin{remark}
  A basis of the space of harmonic polynomials of degree $\ell+1$ is known in
  closed form and is given by the recurrence relation (see
  \cite{PerotChartrand2021}). Notice that the requirement of zero
  integral in $\HPoly{\ell+1}{E}$ can be disregarded in practice, since
  enforcing zero integral into basis functions would not change the results of
  the required computations.
\end{remark}


\section{Numerical Results}

In this section, we propose some numerical experiments to validate our method. We first give numerical evidence of the coercivity of our local bilinear form, 
then we present some convergence tests that asses the theoretical estimates and compare the errors
  {\small
\begin{equation}
    e_{0} \!=\! \frac{\sqrt{\displaystyle\sum_{E \in \mathcal{T}_h} \norm[E]{ u - \proj[\nabla,E]{1}{} \!u_h}^2}}{\norm[\Omega]{u}},
    \quad
    e_{1} \!=\! \frac{\sqrt{\displaystyle\sum_{E \in \mathcal{T}_h} \norm[E]{\sqrt{\K} \left(\nabla u - \nabla \proj[\nabla,E]{1}{} \!u_h\right)^2}}}{\norm[\Omega]{\sqrt{\K}\nabla u}} \,,
    \label{eq:errorsRel}
  \end{equation}
  }
  with respect to the one made by the standard Virtual Element Method \cite{Beirao2015}.


\begin{table}[ht]
\resizebox{\textwidth}{!}{
  \begin{tabular}
      {cccccc}
      \includegraphics[width=1in]{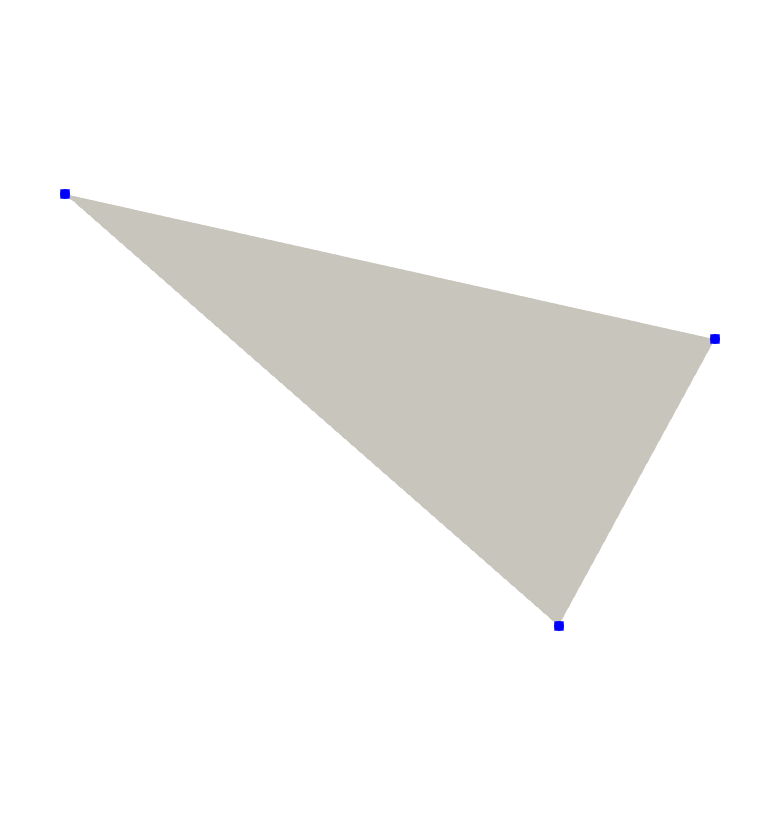} & 
      \includegraphics[width=1in]{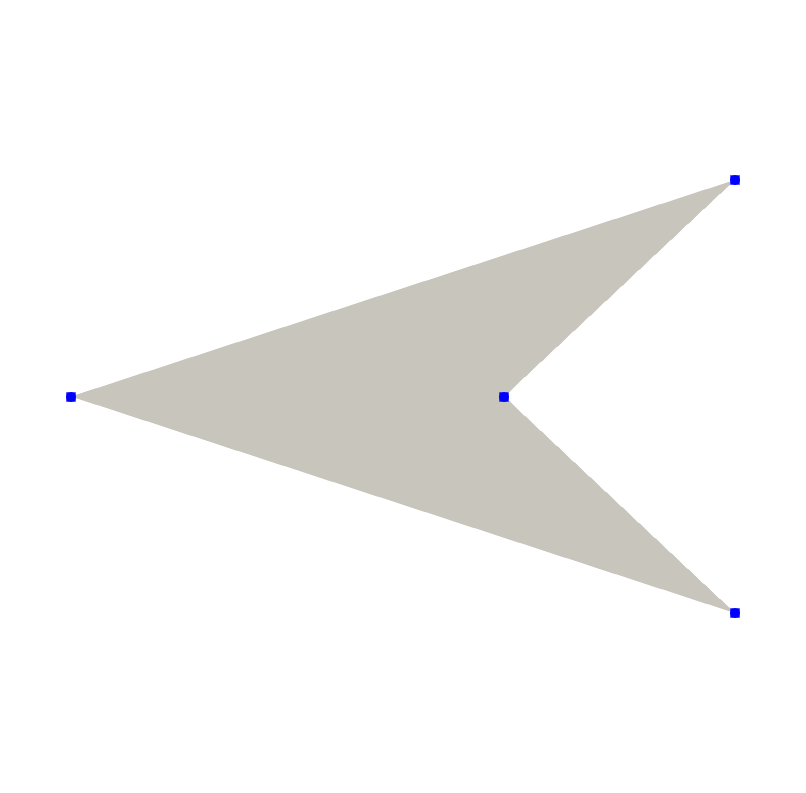} & 
      \includegraphics[width=1in]{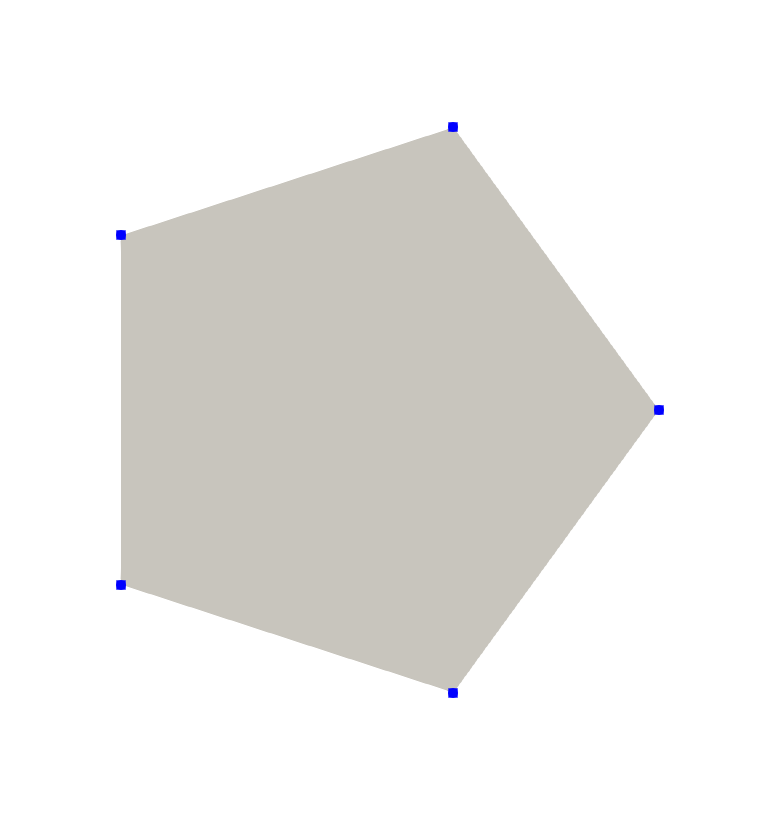} &
      \includegraphics[width=1in]{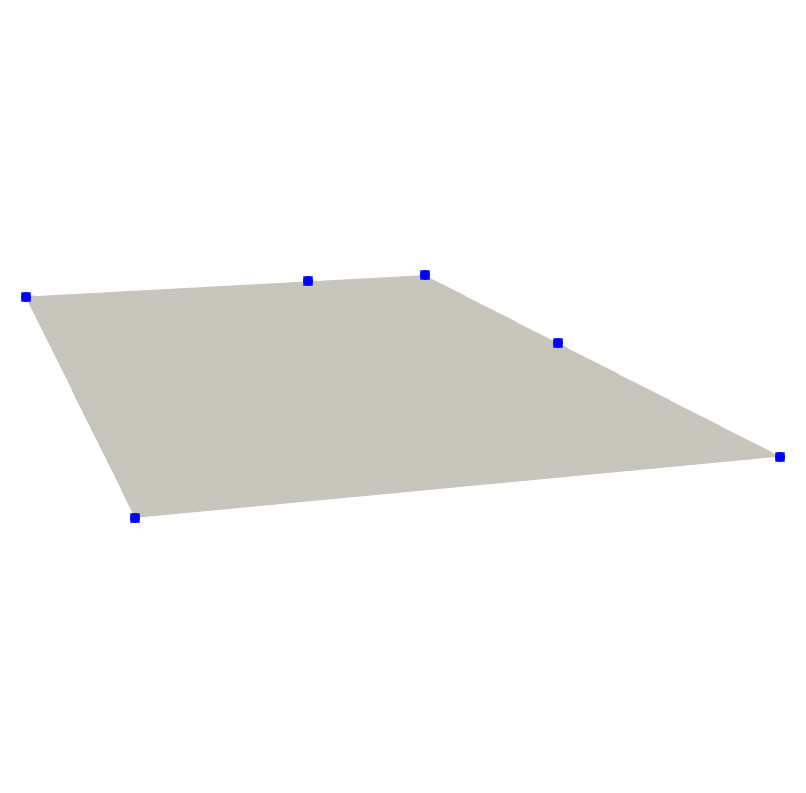} &
      \includegraphics[width=1in]{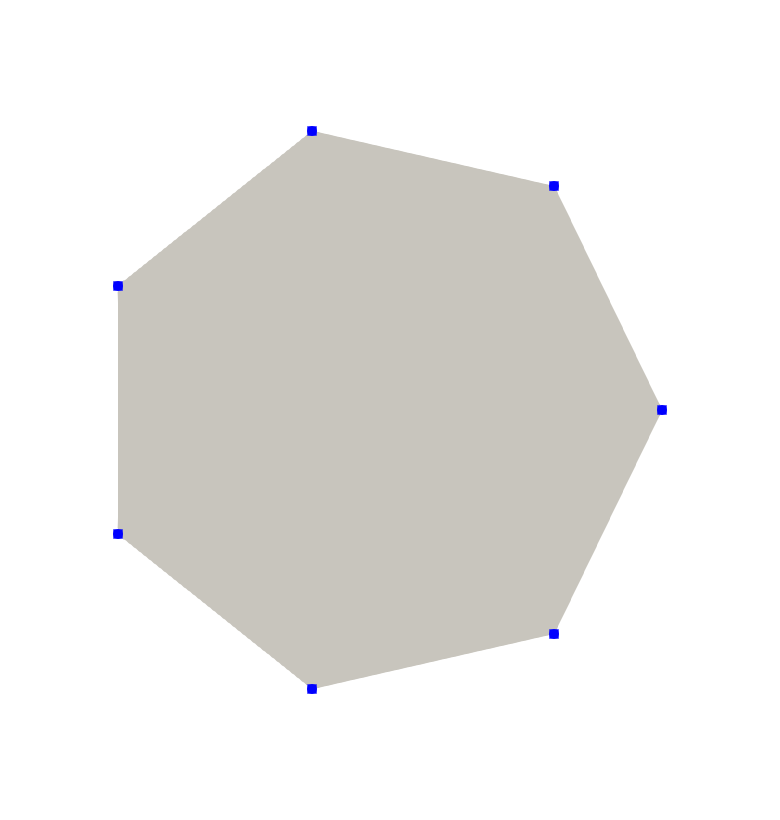} &
      \includegraphics[width=1in]{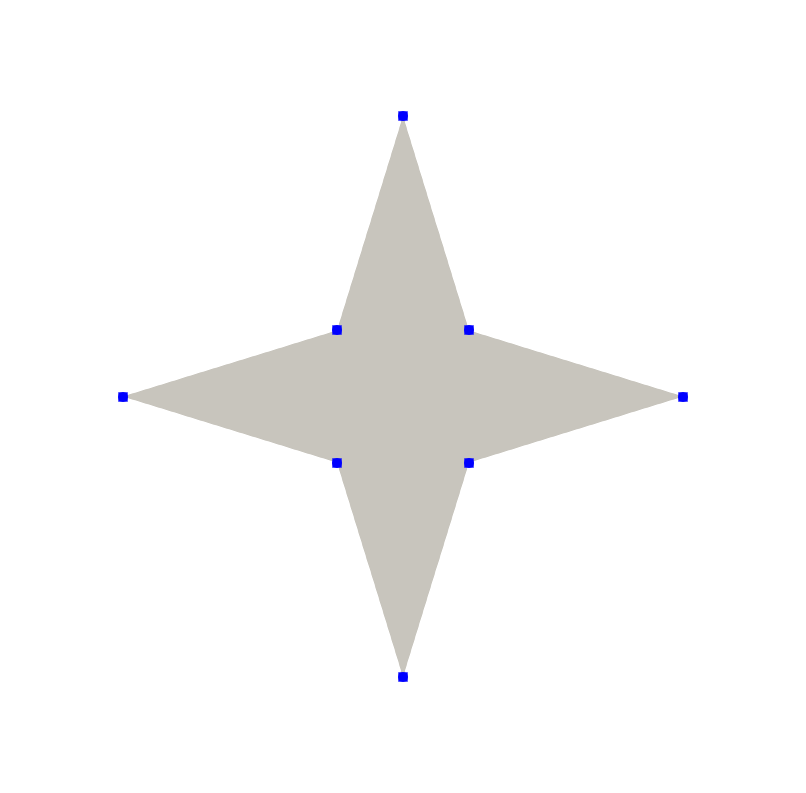}\\
      \multicolumn{1}{c}{Irregular} & \multicolumn{1}{c}{Concave} & \multicolumn{1}{c}{Regular} & \multicolumn{1}{c}{\begin{tabular}[c]{@{}c@{}}Irregular\\ with hanging nodes\end{tabular}} & \multicolumn{1}{c}{Regular} & \multicolumn{1}{c}{Star} \\
      $N_E = 3$, $\ell_E = 0$ & $N_E = 4$, $\ell_E = 1$ & $N_E = 5$, $\ell_E = 1$ & $N_E = 6$, $\ell_E = 2$ & $N_E = 7$, $\ell_E = 2$ & $N_E = 8$, $\ell_E = 3$ \\
      $\sigma_r = 3.8227e\text{-}01
$ & $\sigma_r =1.9207\text{-}01
$ & $\sigma_r =7.1889e\text{-}01$ & $\sigma_r =1.6542e\text{-}01$ & $\sigma_r =6.6611e\text{-}01$ & $\sigma_r = 2.0525e\text{-}01
$\\
      \includegraphics[width=1in]{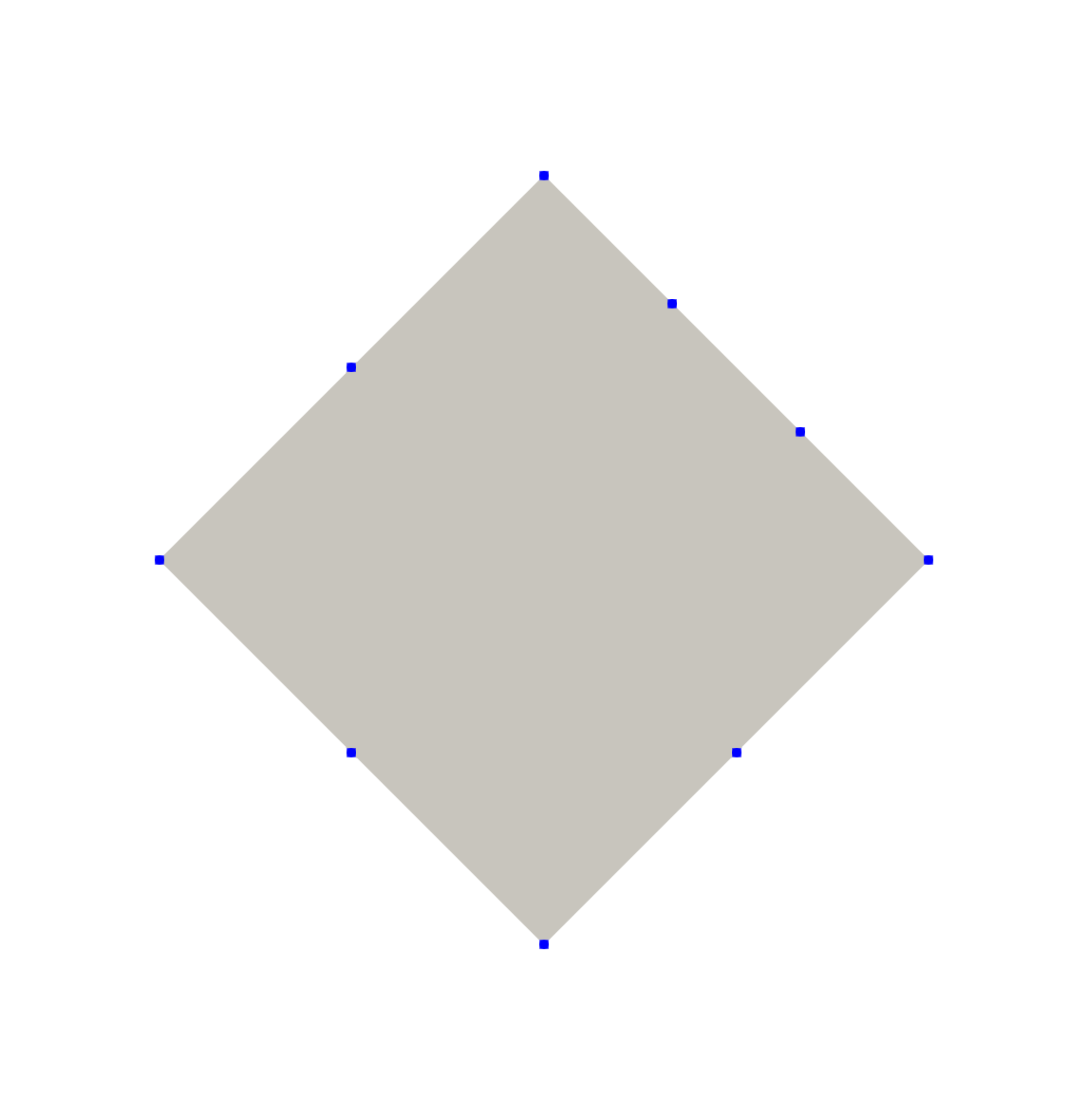} & 
      \includegraphics[width=1in]{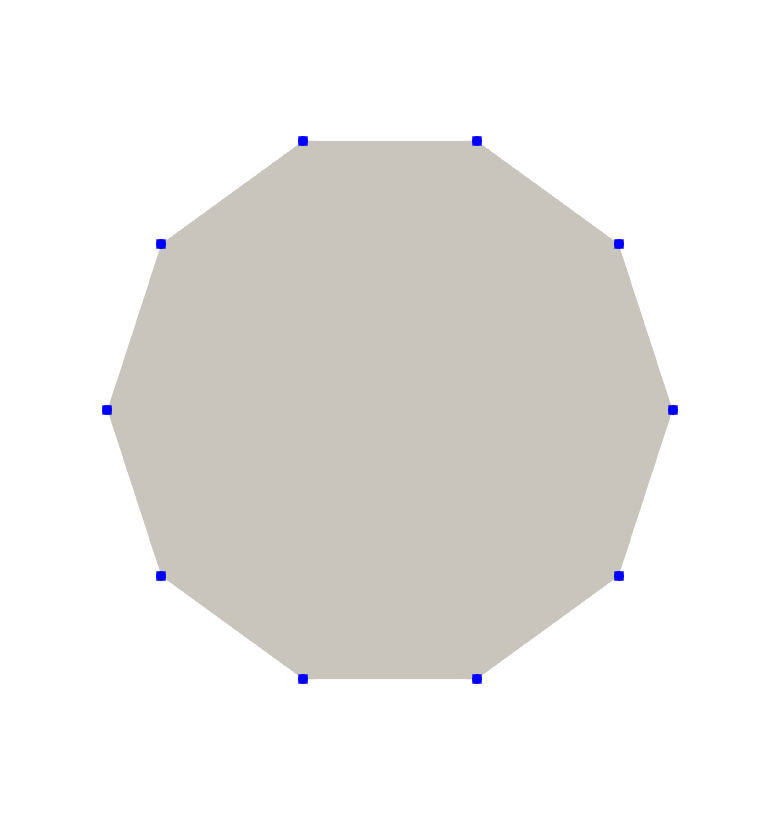} & 
      \includegraphics[width=1in]{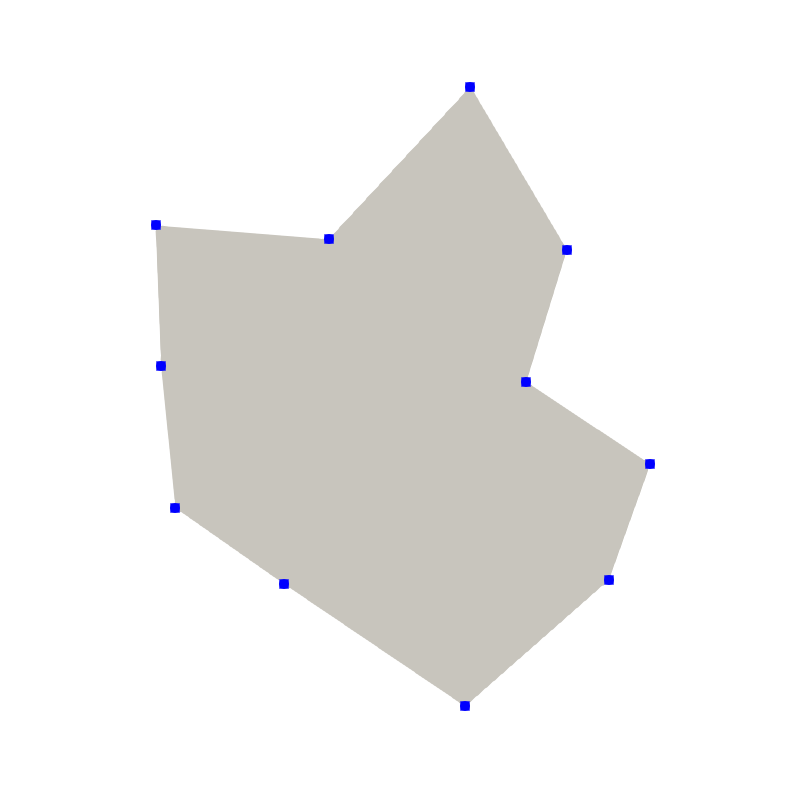} &
      \includegraphics[width=1in]{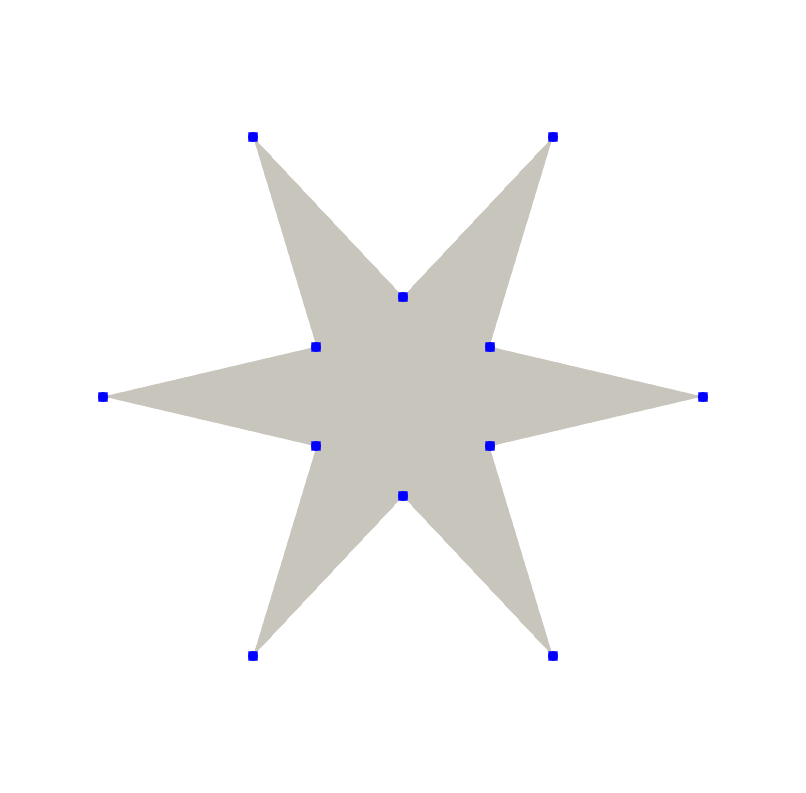} &
      \includegraphics[width=1in]{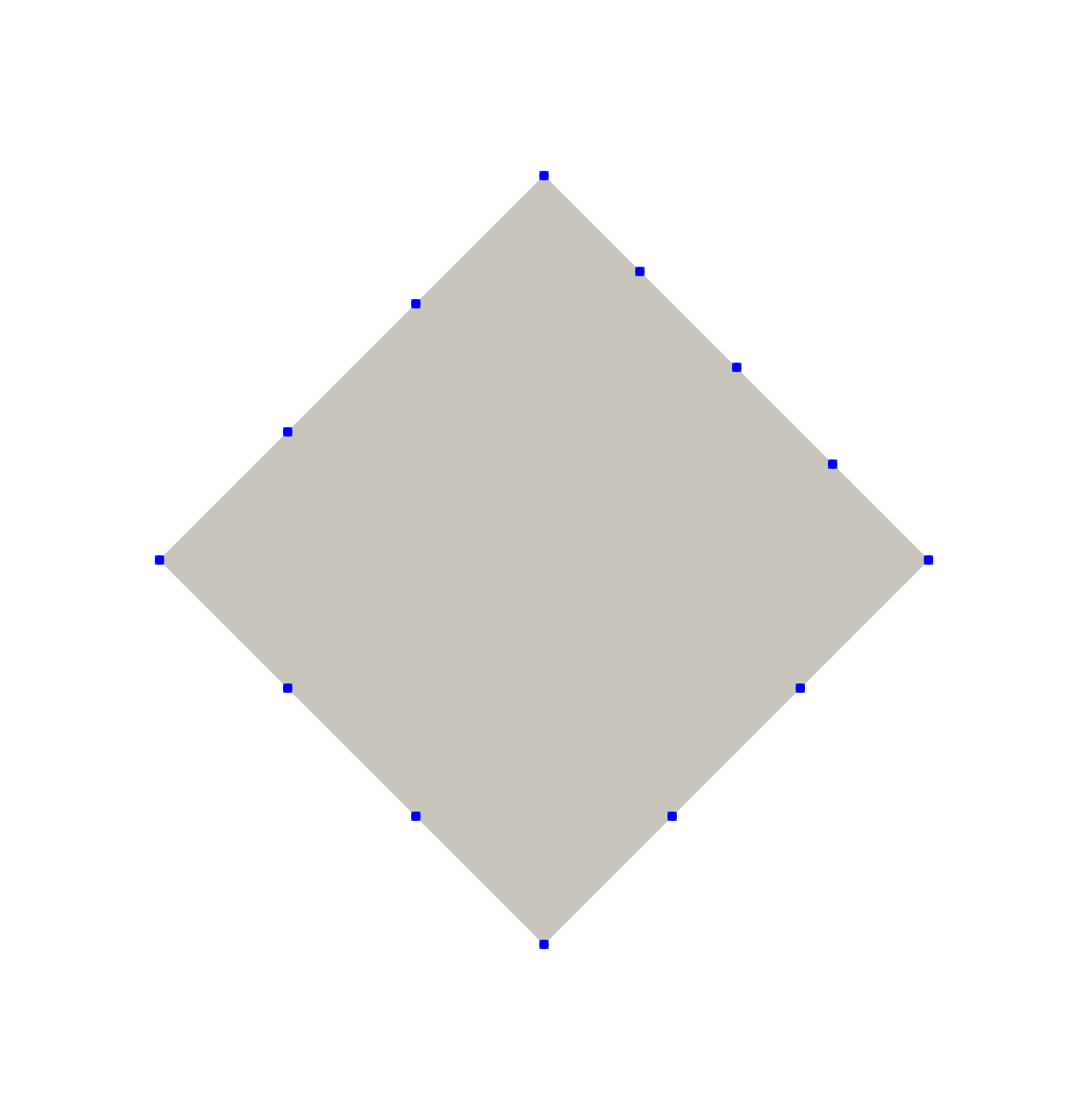} &
      \includegraphics[width=1in]{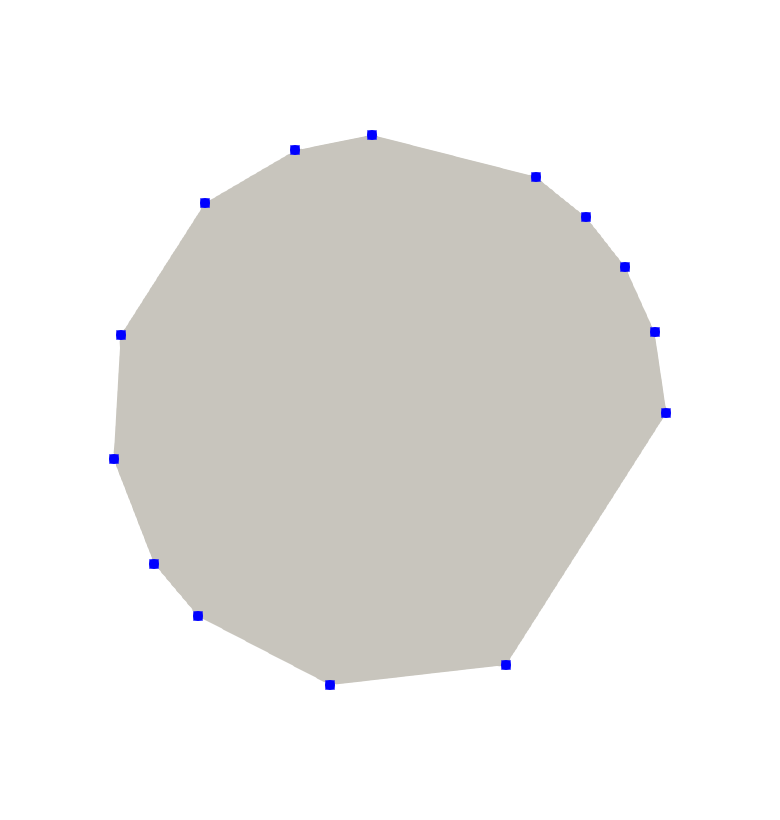}\\
      \multicolumn{1}{c}{\begin{tabular}[c]{@{}c@{}}Irregular\\ with hanging nodes\end{tabular}} & \multicolumn{1}{c}{Regular} & \multicolumn{1}{c}{Concave} & \multicolumn{1}{c}{Star} & \multicolumn{1}{c}{\begin{tabular}[c]{@{}c@{}}Irregular\\ with hanging nodes\end{tabular}} & \multicolumn{1}{c}{Irregular} \\
      $N_E = 9$, $\ell_E = 3$ & $N_E = 10$, $\ell_E = 4$ & $N_E = 11$, $\ell_E = 4$ & $N_E = 12$, $\ell_E = 5$ & $N_E = 13$, $\ell_E = 5$ & $N_E = 14$, $\ell_E = 6$ \\
      $\sigma_r = 2.4452e\text{-}01$ & $\sigma_r =5.8778e\text{-}01$ & $\sigma_r = 1.1917e\text{-}01$ & $\sigma_r =1.0911e\text{-}01$ & $\sigma_r =1.5378e\text{-}01$ & $\sigma_r = 4.8291e\text{-}02
$\\
      \includegraphics[width=1in]{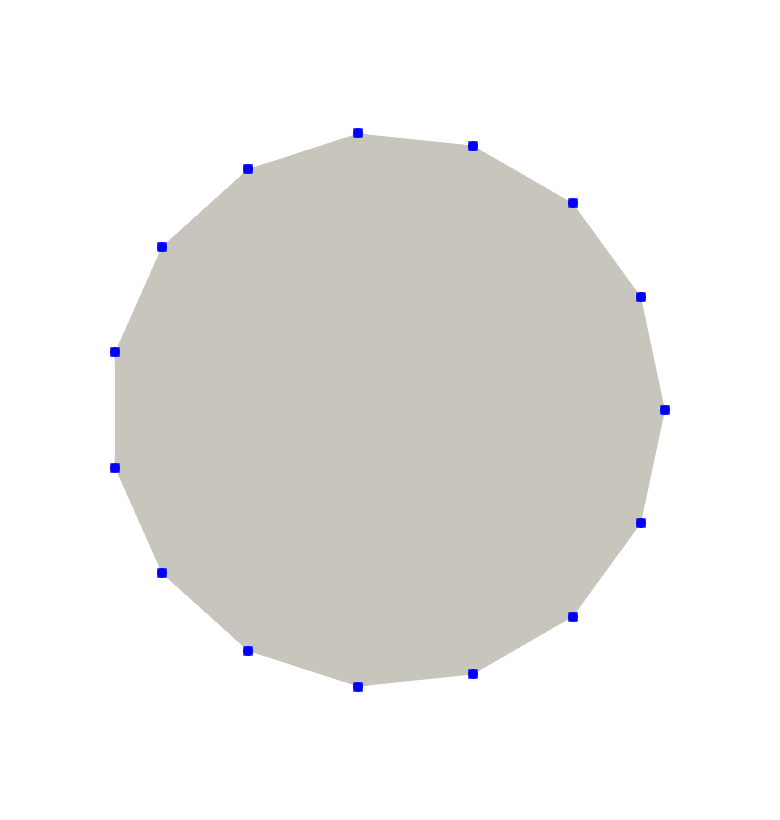} & 
      \includegraphics[width=1in]{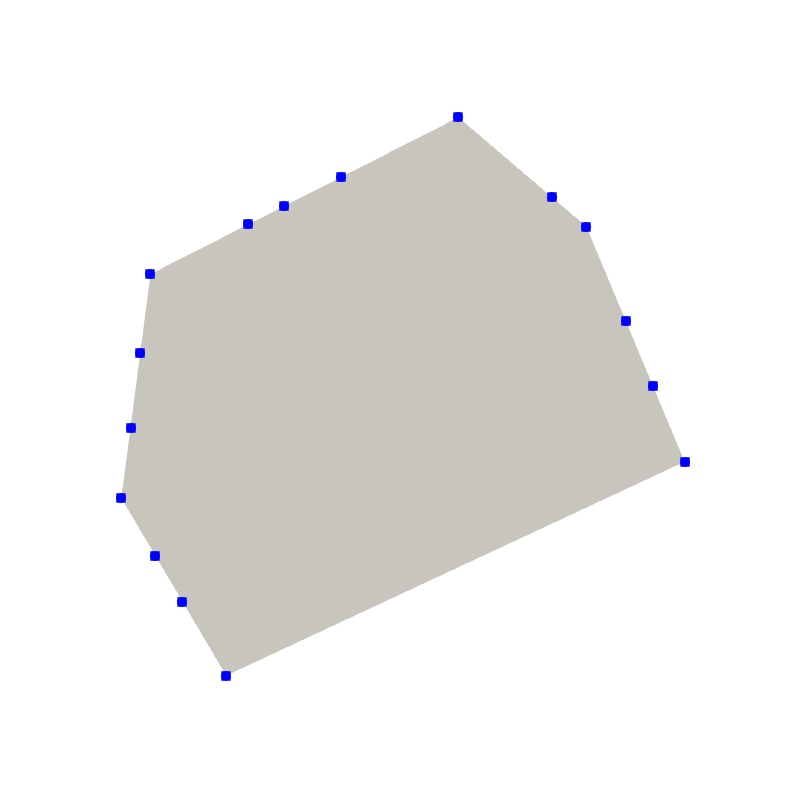} & 
      \includegraphics[width=1in]{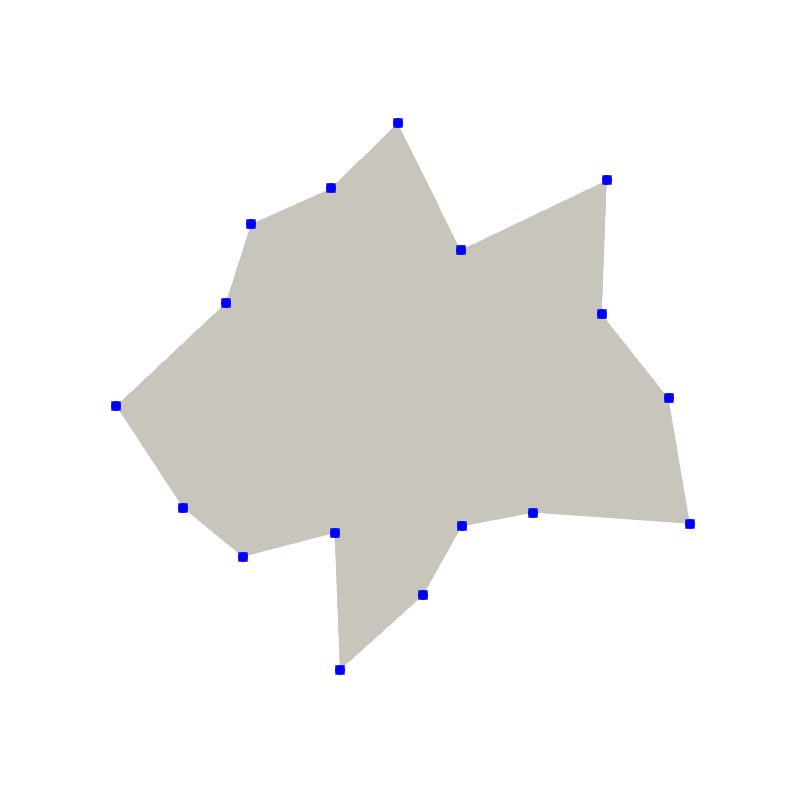} &
      \includegraphics[width=1in]{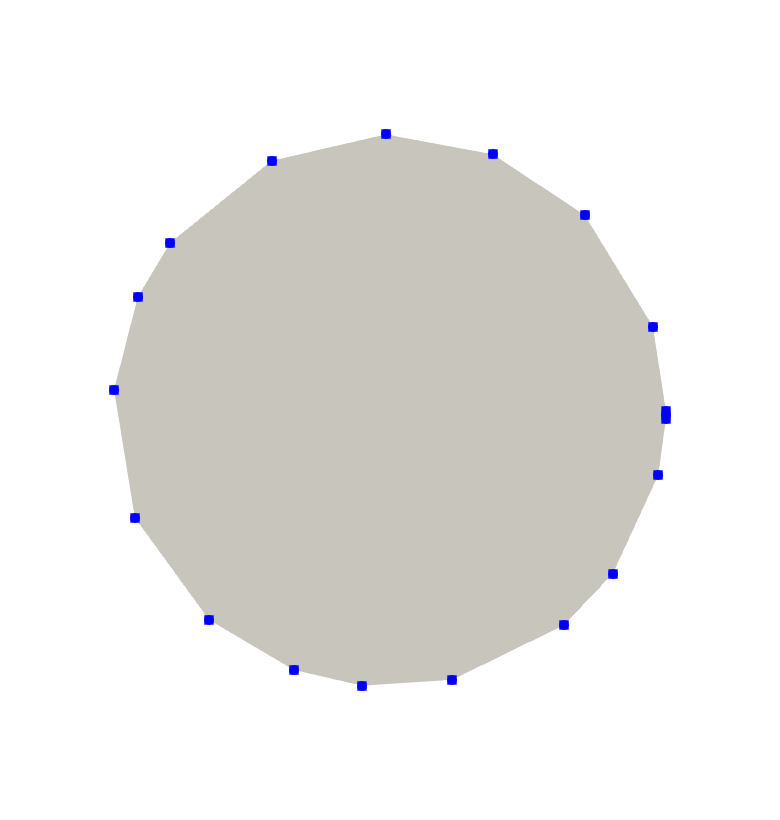} &
      \includegraphics[width=1in]{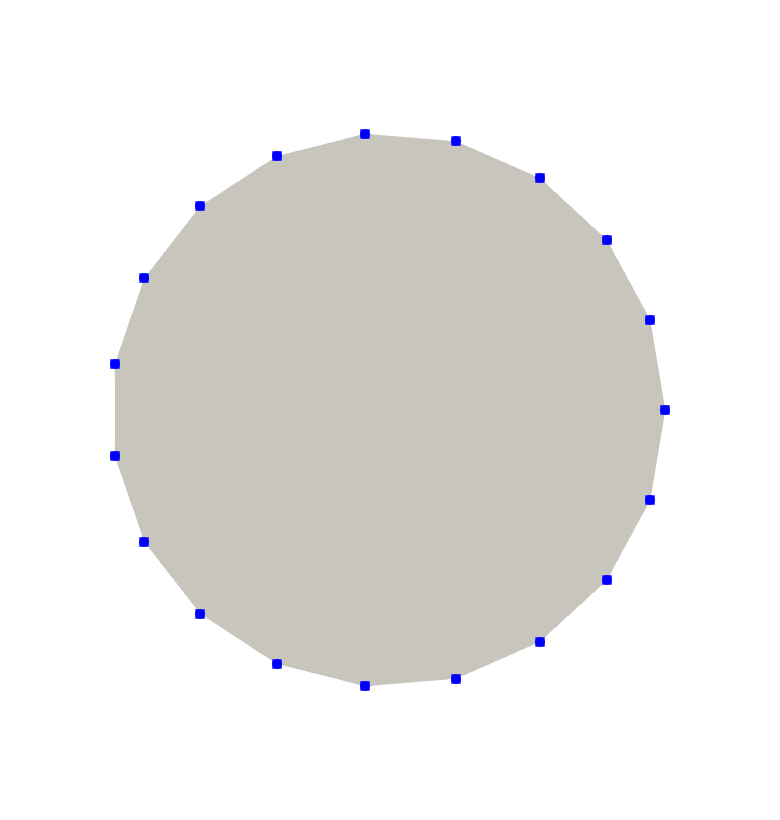} &
      \includegraphics[width=1in]{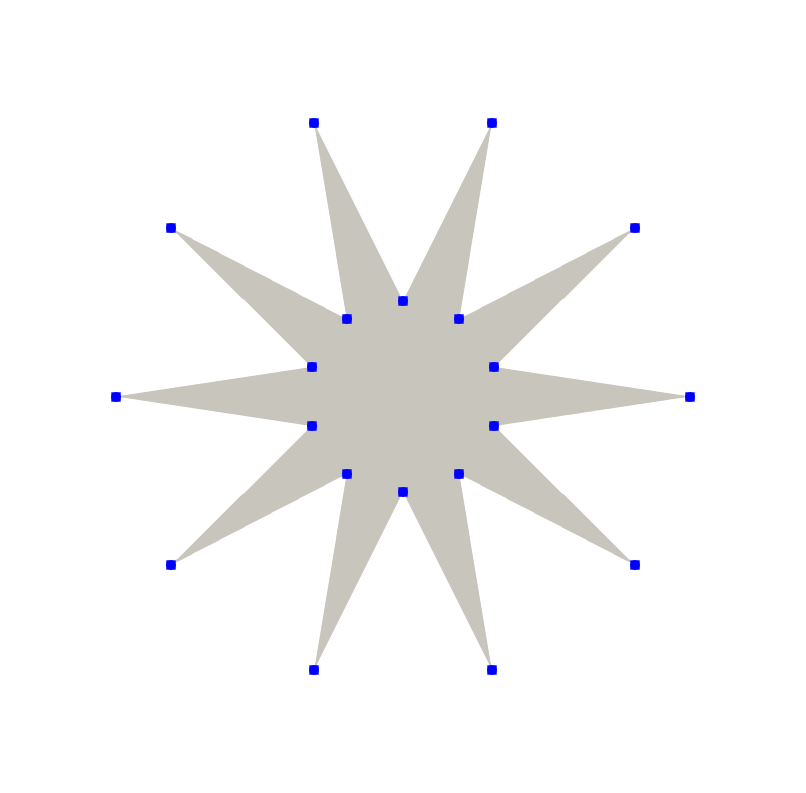}\\
      \multicolumn{1}{c}{Regular} & \multicolumn{1}{c}{\begin{tabular}[c]{@{}c@{}}Irregular\\ with hanging nodes\end{tabular}}  & \multicolumn{1}{c}{Concave} & \multicolumn{1}{c}{\begin{tabular}[c]{@{}c@{}}Irregular\\ with a collapsing edge\end{tabular}} & \multicolumn{1}{c}{Regular} & \multicolumn{1}{c}{Star} \\
      $N_E = 15$, $\ell_E = 6$ & $N_E = 16$, $\ell_E = 7$ & $N_E = 17$, $\ell_E = 7$ & $N_E = 18$, $\ell_E = 8$ & $N_E = 19$, $\ell_E = 9$ & $N_E = 20$, $\ell_E = 10$ \\
      $\sigma_r =4.0674e\text{-}01$ & $\sigma_r = 1.3047e\text{-}04
$ & $\sigma_r =1.1031e\text{-}02$ & $\sigma_r = 2.3334e\text{-}02
$ & $\sigma_r =3.2470e\text{-}01$ & $\sigma_r =3.6314e\text{-}02$\\
  \end{tabular}
  }
  \caption{$\sigma_r$ of the elemental stiffness matrices related to different kinds of polygons.} 
  \label{tab:singularvalues}
\end{table}

In the first test, we consider a set of different polygons, with different geometrical features, such as concavities, symmetries,  and aligned edges. For each polygon, choosing $\ell_E$ according to Theorem \ref{th:wellposed}, we asses the local stability of the discrete diffusion operator \eqref{eq:discrvarform} ($\K=1$, $\bs\beta=0$, and $\gamma=0$), evaluating the second smallest singular value of the stiffness matrix  denoted by $\sigma_r$.
The results, reported in Table \ref{tab:singularvalues}, confirm the stability of the method and good robustness with respect to the geometrical complexity being $\sigma_r$ always well detached from zero (the smallest singular value of the stiffness matrix is always vanishing).

  \begin{figure}
    \centering	
    \subfigure[\label{fig:Mesh_QD_60x60}]	{\includegraphics[width=.32\textwidth, height = .17\textheight]{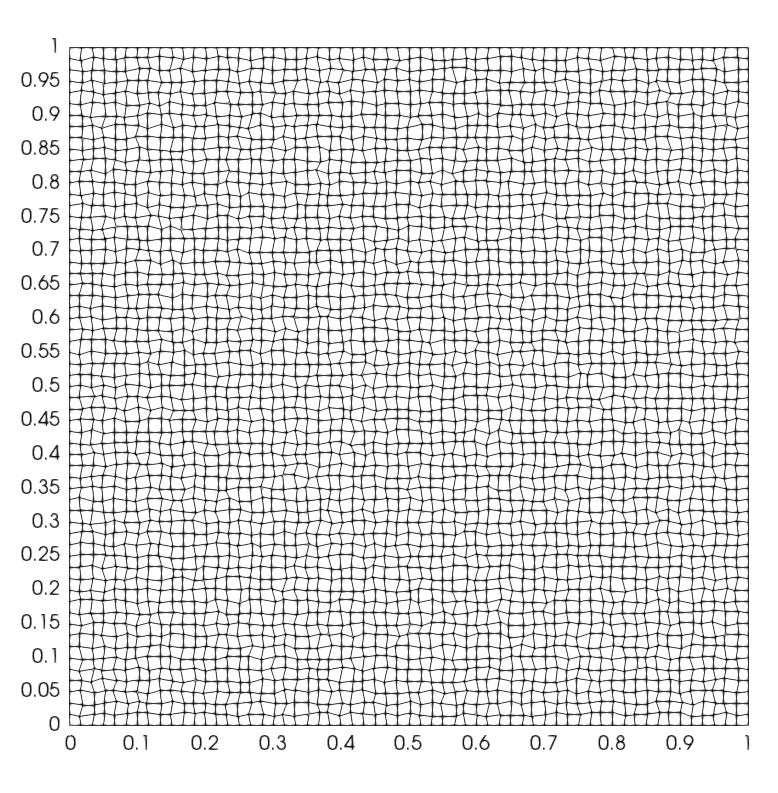}}	\subfigure[\label{fig:Mesh_ED_44x50}]	{\includegraphics[width=.32\textwidth, height = .17\textheight]{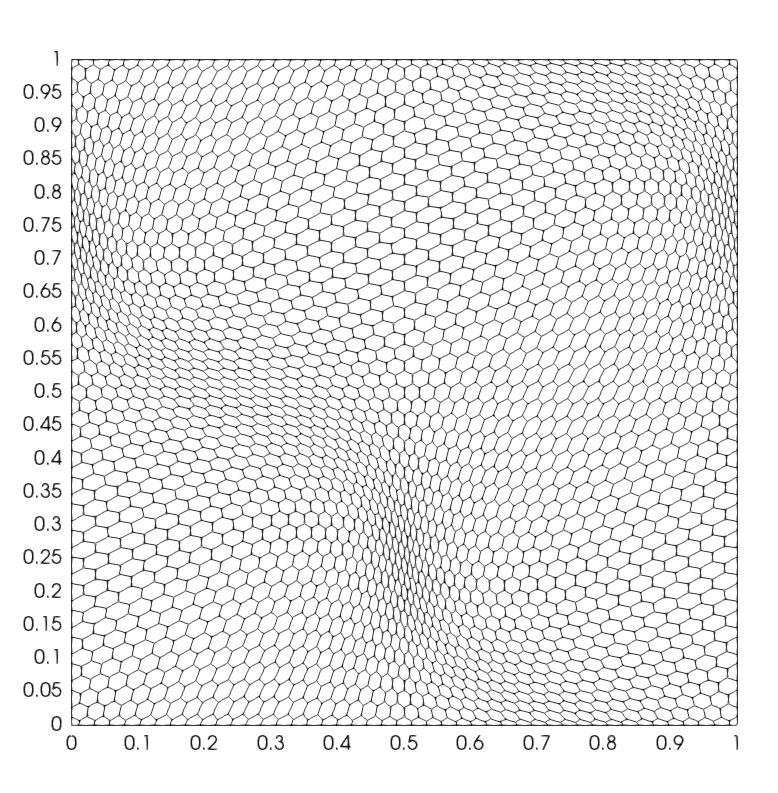}}	\subfigure[\label{fig:Mesh_EDD_44x50}]	{\includegraphics[width=.32\textwidth, height = .17\textheight]{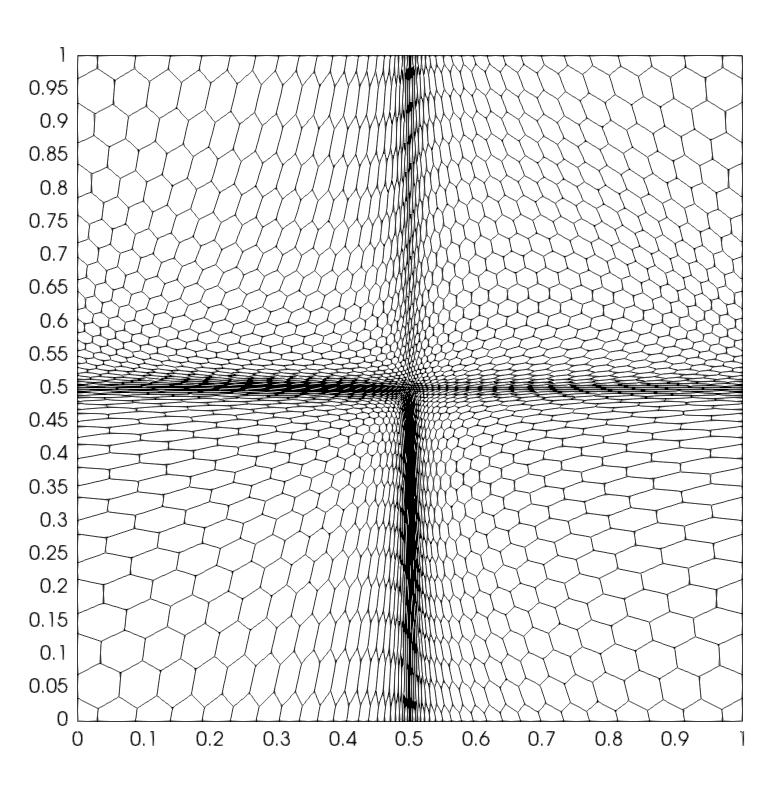}}
    \caption{Meshes used in the numerical experiments. Left: Distorted squared mesh. Center: Distorted Voronoi mesh. Right: Highly-distorted Voronoi mesh.}
    \label{fig:mesh}
  \end{figure}

In the second test, we compare the stabilization-free Virtual Element Method (SFVEM in short) with the standard VEM with the \textit{dofi-dofi} stabilization term (VEM in short) \cite{Beirao2015b} by plotting the relative errors $e_0$ and $e_1$ \eqref{eq:errorsRel}, and computing their rates of convergence on 
three families of distorted and highly-distorted meshes.
The fourth refinement of each family of meshes is shown in Figure \ref{fig:mesh}.
  In order to show the advantages of SFVEM with respect to the standard VEM, as suggested in \cite{BERRONE2022107971}, we consider an \textit{anisotropic} diffusion tensor $\K$. Let $\Omega$ be the unit square, we consider the advection-diffusion-reaction problem \eqref{eq:modelvarform} with coefficients
  
%
\begin{equation*}
\K = \GG \begin{bmatrix}
1 & 0\\
0 & 1.0e\text{-}09
\end{bmatrix} \GG^T,
\quad\GG = \begin{bmatrix}
\cos(\theta) & -\sin(\theta)\\
\sin(\theta) & \cos(\theta)
\end{bmatrix},
\quad \bbeta(x,y) = \begin{bmatrix}
\bbeta_1(x,y)\\
\bbeta_2(x,y)
\end{bmatrix},
\end{equation*}
and $\gamma(x,y) = x(1-x)y(1-y)$,
where $\GG$ is the Givens rotation matrix with $\theta \in \mathbb{R}$.
For $R_1$, $R_2 \in [0,1]$, we define \cite{SBEnumath}
\begin{align*}
  \bbeta_1(x,y;R_1)&=
                 250000 x^4 y^3 (R_1-x) (1-x)^4
  \\ \nonumber
               &\left[4 R_2 \left(1-5 y
                 +9 y^2 -7 y^3 +2 y^4\right)
                 -5 y+24 y^2-42 y^3
                 +32 y^4-9 y^5\right],\\ \nonumber
  \bbeta_2(x,y;R_2)&=-\bbeta_1(y,x;R_2),
\end{align*}
and we fix $R_1 = 0.9$, $R_2 = 0.3$ and $\theta = \frac{\pi}{6}$. We choose $f(x,y)$ in such a way the exact solution is $u(x,y)= \bbeta_1(x,y)$.

  \begin{figure}
    \centering
    \subfigure[\label{fig:ErrorL2_QD_EP_R1:0.900000_R2:0.300000_Theta:0.523560}]
    {\includegraphics[width=.32\textwidth, height = .17\textheight]{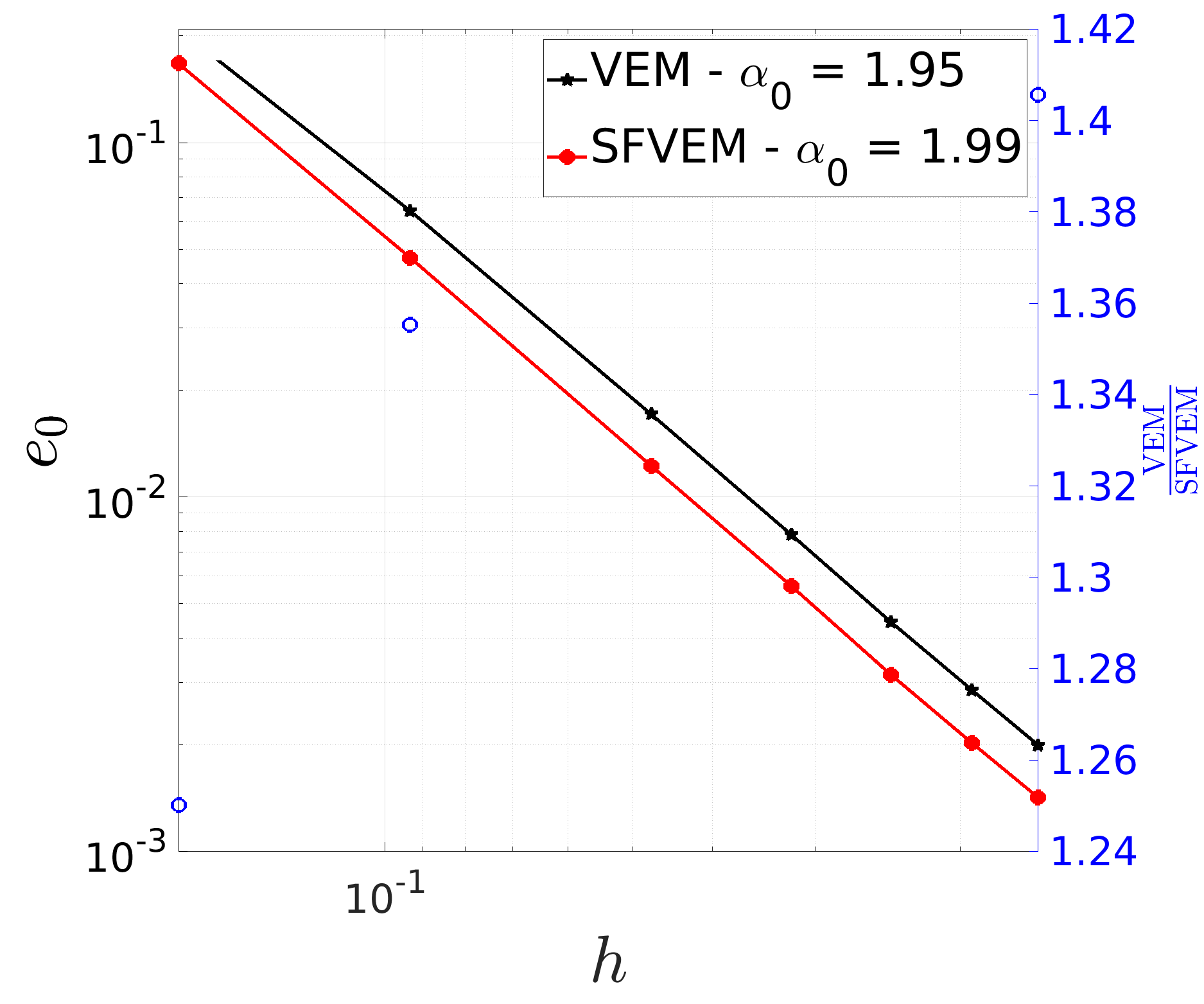}}	\subfigure[\label{fig:ErrorL2_ED_EP_R1:0.900000_R2:0.300000_Theta:0.523560}]
    {\includegraphics[width=.32\textwidth, height = .17\textheight]{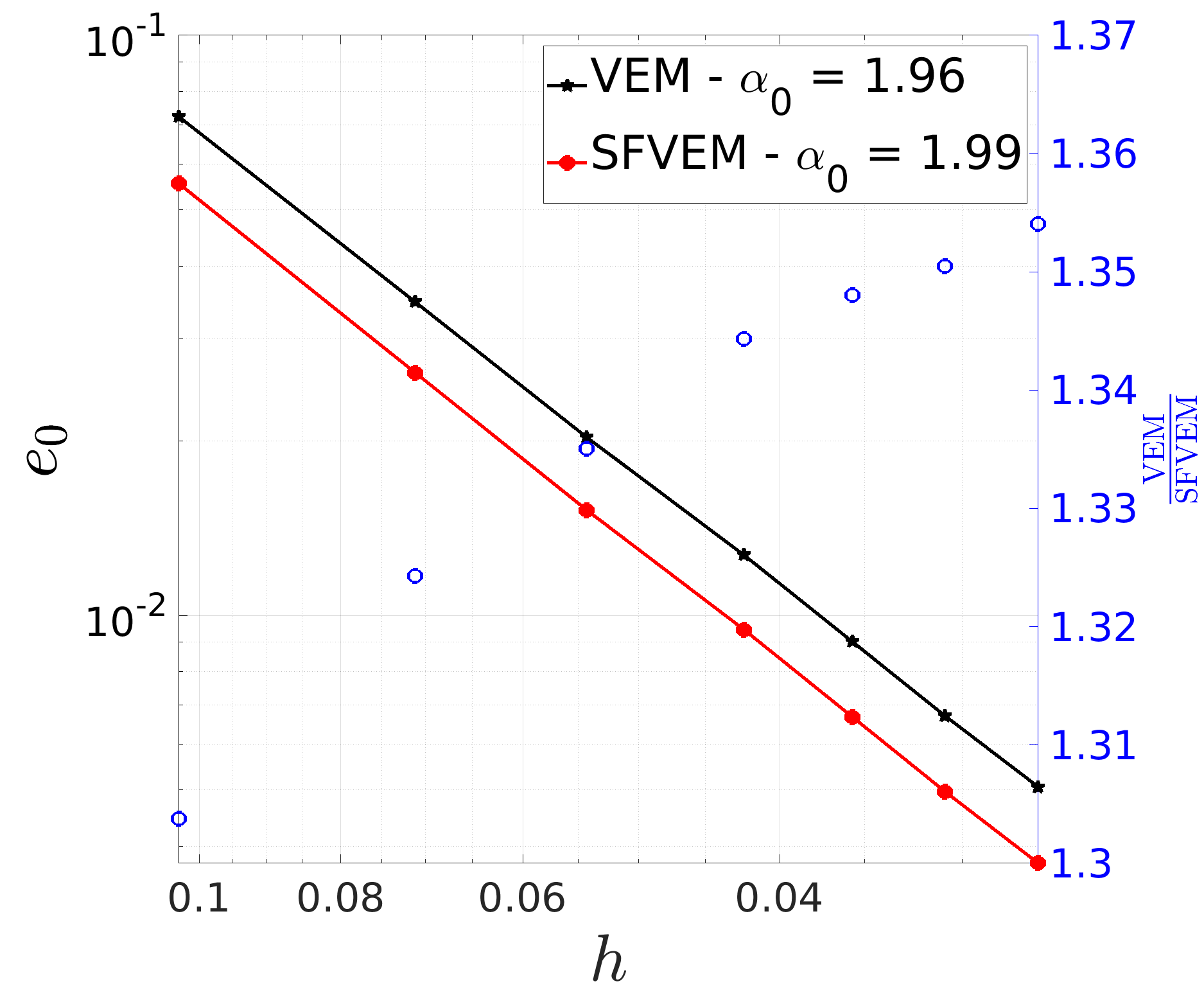}}
    \subfigure[\label{fig:ErrorL2_EDD_EP_R1:0.900000_R2:0.300000_Theta:0.523560}]
    {\includegraphics[width=.32\textwidth, height = .17\textheight]{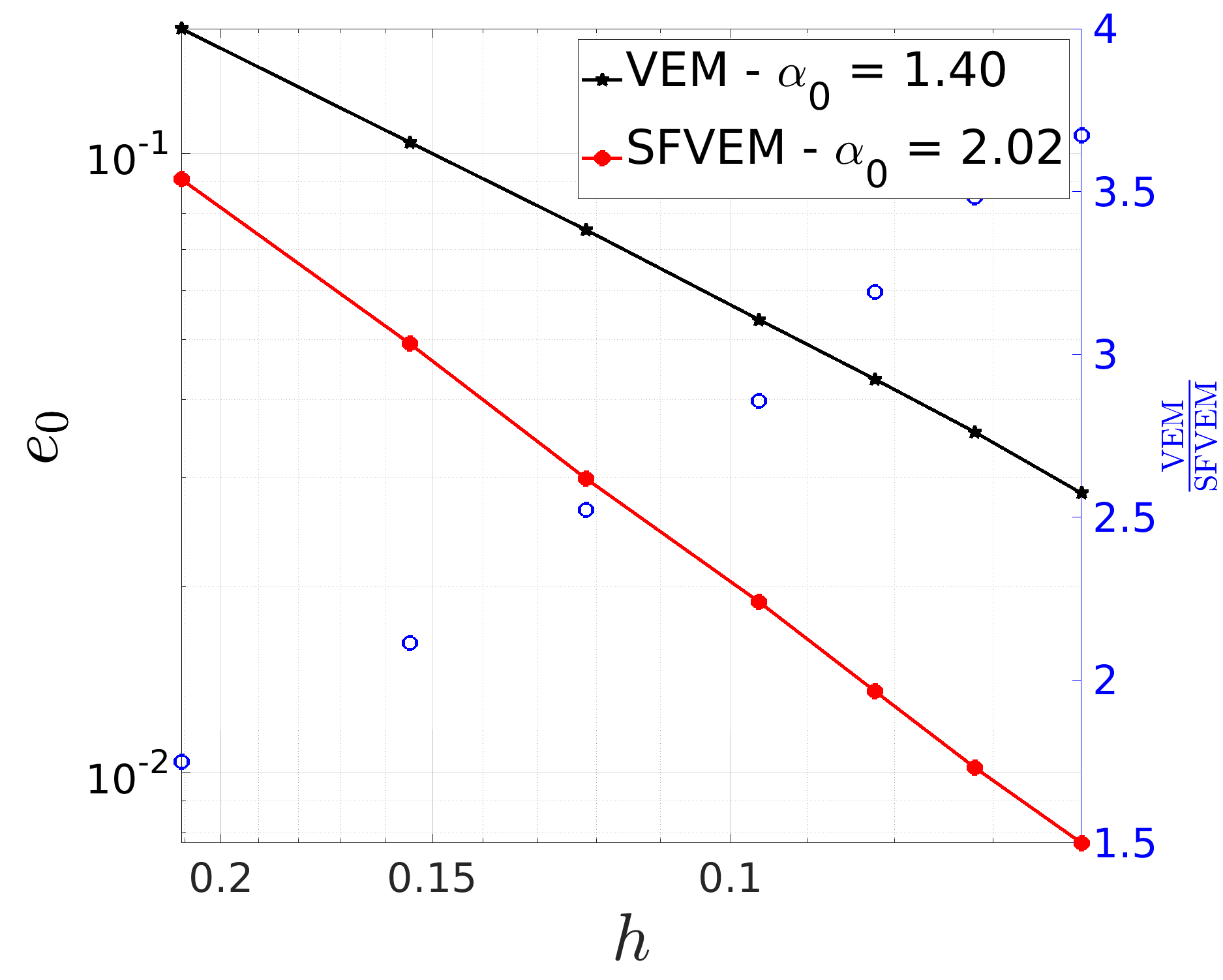}}	\subfigure[\label{fig:ErrorH1_QD_EP_R1:0.900000_R2:0.300000_Theta:0.523560}]
    {\includegraphics[width=.32\textwidth, height = .17\textheight]{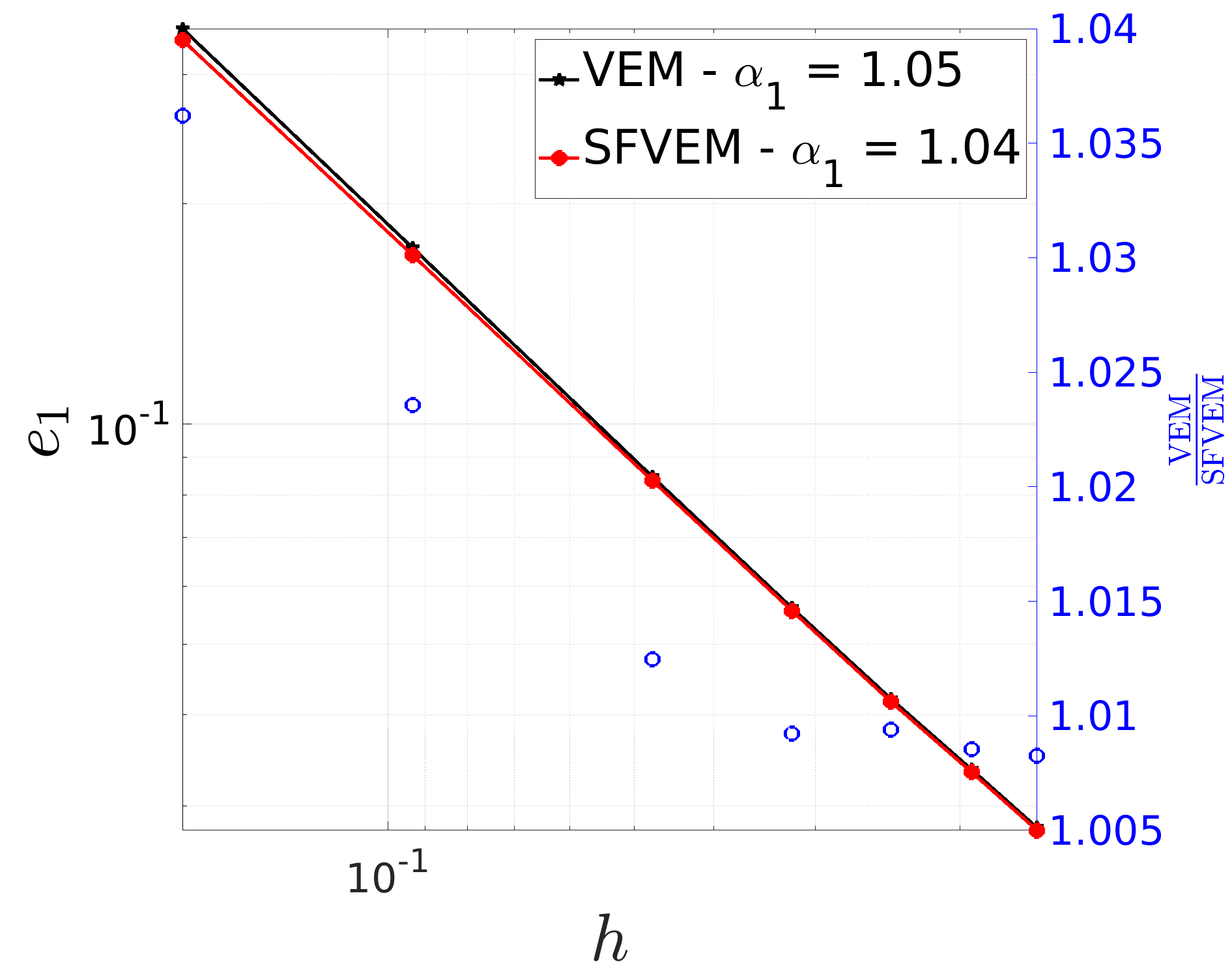}}	\subfigure[\label{fig:ErrorH1_ED_EP_R1:0.900000_R2:0.300000_Theta:1.570796}]
    {\includegraphics[width=.32\textwidth, height = .17\textheight]{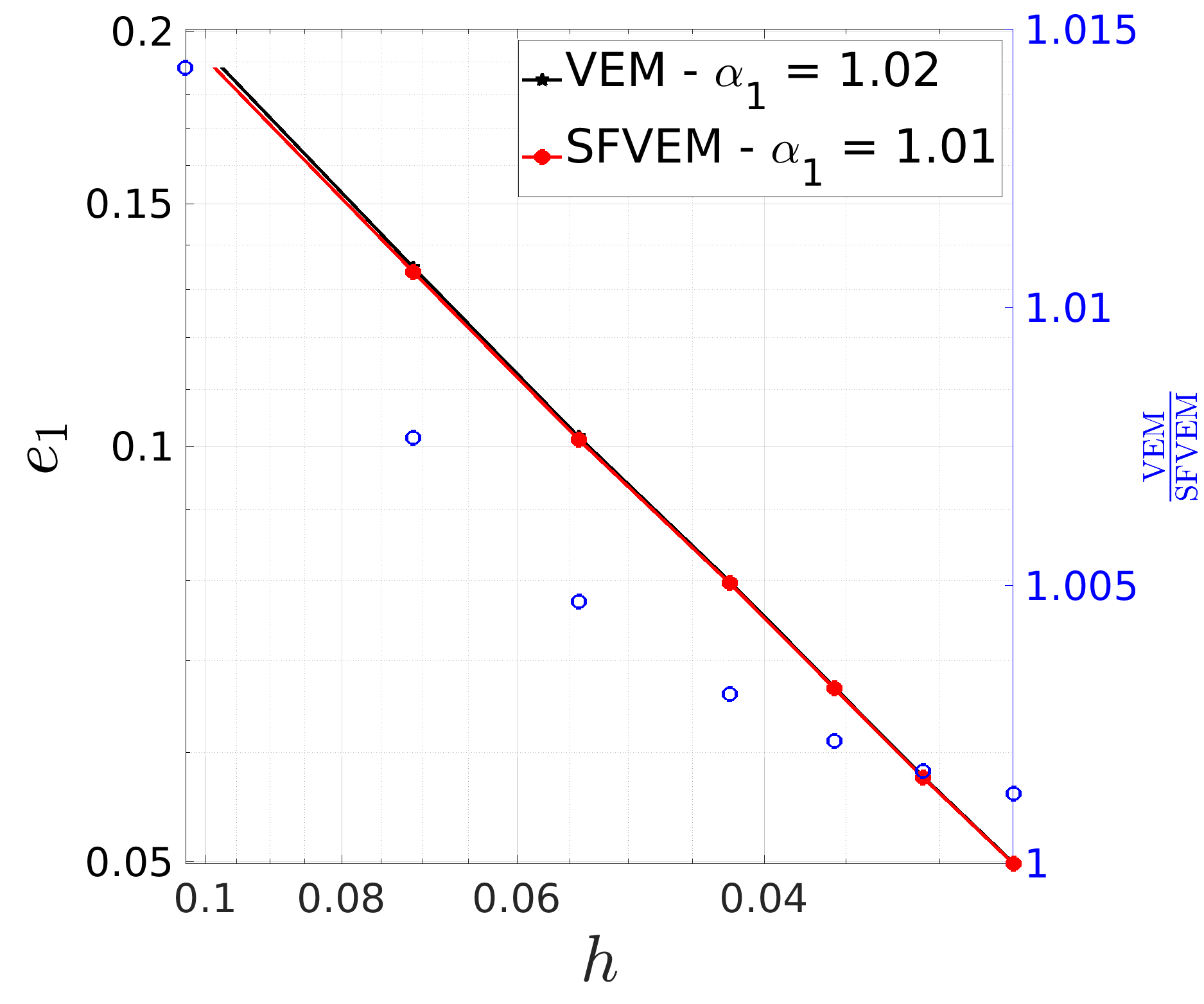}}	\subfigure[\label{fig:ErrorH1_EDD_EP_R1:0.900000_R2:0.300000_Theta:1.570796}]
    {\includegraphics[width=.32\textwidth, height = .17\textheight]{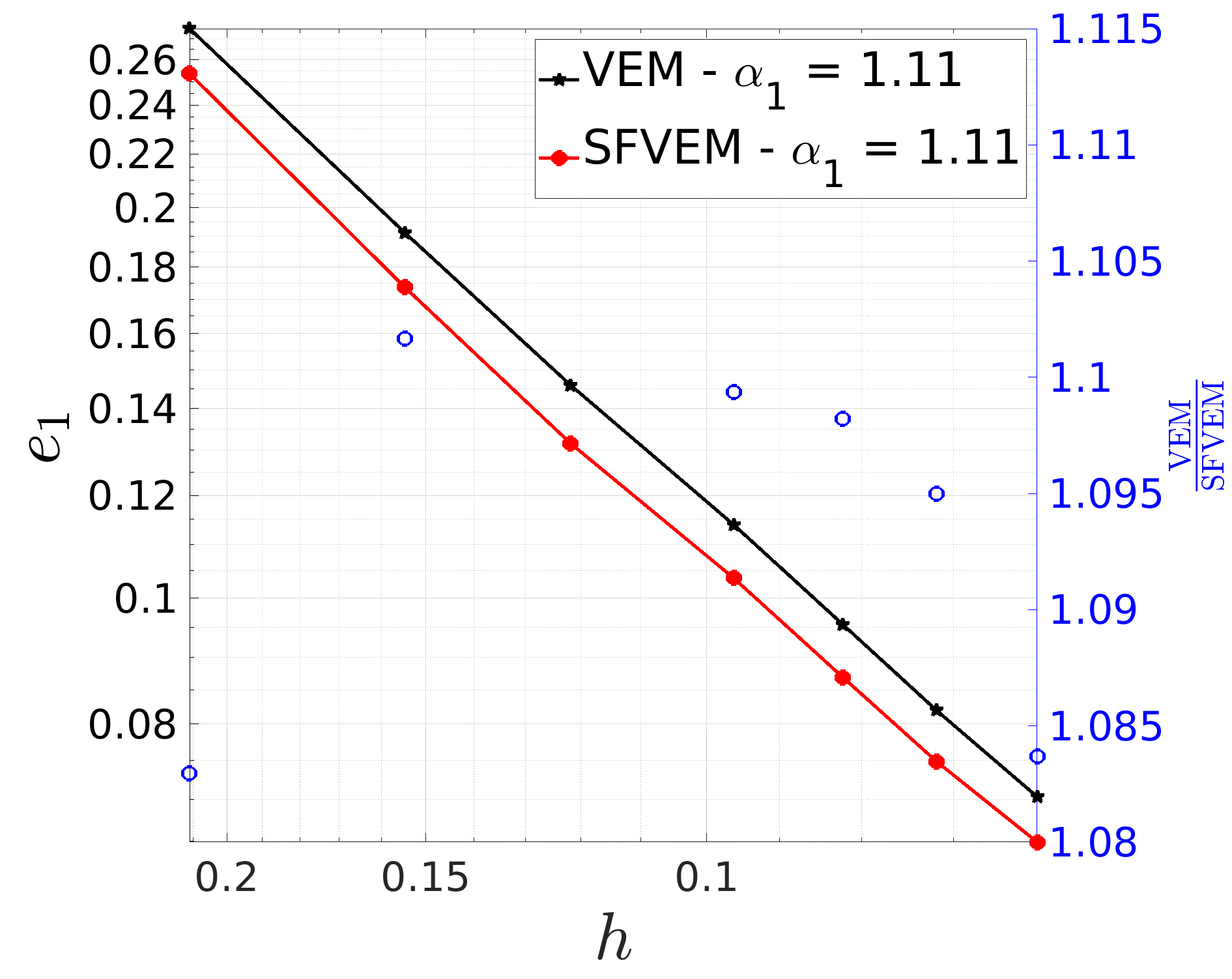}}
    \caption{Behaviour of errors $e_0$ and $e_1$ \eqref{eq:errorsRel} w.r.t. $h$. Left: Distorted squared mesh. Center: Distorted Voronoi mesh. Right: Highly-distorted Voronoi mesh.}
    \label{fig:errors}
  \end{figure}

  In Figure \ref{fig:errors}, we plot the convergence curves of errors $e_0$ and $e_1$ \eqref{eq:errorsRel} and the ratio between their values for VEM and SFVEM (right axis of each figure). The legends report the rates of convergence of the errors ($\alpha_0$ and $\alpha_1$, respectively). 
  The performances of the two methods are almost equivalent concerning the $e_1$ error, see Figures  \ref{fig:ErrorH1_QD_EP_R1:0.900000_R2:0.300000_Theta:0.523560}-\ref{fig:ErrorH1_EDD_EP_R1:0.900000_R2:0.300000_Theta:1.570796}.
  Whereas in Figures \ref{fig:ErrorL2_QD_EP_R1:0.900000_R2:0.300000_Theta:0.523560}-\ref{fig:ErrorL2_EDD_EP_R1:0.900000_R2:0.300000_Theta:0.523560} SFVEM easily reaches the asymptotic rates of convergence on all the meshes and displays a smaller $e_0$ error, whereas VEM is still in a pre-asymptotic regime on highly-distorted Voronoi meshes  and displays an error between two and three times w.r.t. SFVEM.

\label{sec:numres}

\section{Conclusion}
We propose a new first-order stabilization-free VEM that exploits projections on
harmonic polynomials to build a self-stabilized bilinear form. Numerical results
show good stability of the method and optimal rates of convergence.

\bibliographystyle{plain}
\bibliography{biblio}

\end{document}